\documentclass[11pt,twoside]{article}
\usepackage{amssymb}
\usepackage{hyperref}
\usepackage{float}
\usepackage{amsfonts}
\usepackage[centertags]{amsmath}
\usepackage{amsthm}
\usepackage{newlfont}
\usepackage{array}

\date{}

\begin{document}

\centerline{\Large{\bf Integrated and Differentiated Sequence Spaces}}

\centerline{}

\centerline{\Large{\bf and Weighted Mean}}
\centerline{}

\centerline{\bf {Murat Kiri\c{s}ci}}

\centerline{}

\centerline{Department of Mathematical Education, Hasan Ali Y\"{u}cel Education Faculty,}

\centerline{ Istanbul University, Vefa, 34470, Fatih, Istanbul, Turkey}

\centerline{}

\newtheorem{Theorem}{\quad Theorem}[section]

\newtheorem{Definition}[Theorem]{\quad Definition}

\newtheorem{Corollary}[Theorem]{\quad Corollary}

\newtheorem{Lemma}[Theorem]{\quad Lemma}

\newtheorem{Example}[Theorem]{\quad Example}

\newtheorem*{remark}{Remark}

\begin{abstract}
The purpose of this paper is twofold. Firstly, the new matrix domains are constructed with the new infinite matrices and some properties are investigated.
Furthermore, dual spaces of new matrix domains are computed and matrix transformations are characterized. Secondly, examples between new spaces with classical sequence spaces and sequence spaces which are derived by an infinite matrix are given in the table form.
\end{abstract}

\centerline{}



\section{Introduction}

It is well known that, the $\omega$ denotes the family of all real (or complex)-valued sequences.
$\omega$ is a linear space and each linear subspace of $\omega$ (with the included addition
and scalar multiplication) is called a \emph{sequence space} such as the spaces $c$, $c_{0}$ and
$\ell_{\infty}$, where $c$, $c_{0}$ and $\ell_{\infty}$ denote the set of all convergent
sequences in fields $\mathbb{R}$ or $\mathbb{C}$, the set of all null sequences and the set
of all bounded sequences, respectively. It is clear that the sets $c$, $c_{0}$ and $\ell_{\infty}$
are the subspaces of the $\omega$. Thus, $c$, $c_{0}$ and $\ell_{\infty}$ equipped with a vector space structure,
from a sequence space. By $bs$ and $cs$, we define the spaces of all bounded and convergent series, respectively.\\

\emph{A coordinate space} (or \emph{$K-$space}) is a vector space of numerical sequences, where addition and scalar multiplication are defined pointwise. That is, a sequence space $X$ with a linear topology is called a $K$-space provided each of the maps $p_{i}:X\rightarrow \mathbb{C}$ defined by $p_{i}(x)=x_{i}$ is continuous for all $i\in \mathbb{N}$. A $K-$space is called an \emph{$FK-$space} provided $X$ is a complete linear metric space. An \emph{$FK-$space} whose topology is normable is called a \emph{$BK-$ space}. If a normed sequence space $X$ contains a sequence $(b_{n})$ with the property that for every $x\in X$ there is unique sequence of scalars $(\alpha_{n})$ such that
\begin{eqnarray*}
\lim_{n\rightarrow\infty}\|x-(\alpha_{0}b_{0}+\alpha_{1}b_{1}+...+\alpha_{n}b_{n})\|=0
\end{eqnarray*}
then $(b_{n})$ is called \emph{Schauder basis} for $X$. The series $\sum\alpha_{k}b_{k}$ which has the sum $x$ is then called the expansion of $x$ with respect to $(b_{n})$, and written as $x=\sum\alpha_{k}b_{k}$. An \emph{$FK-$space} $X$ is said to have $AK$ property, if $\phi \subset X$ and $\{e^{k}\}$ is a basis for $X$, where $e^{k}$ is a sequence whose only non-zero term is a $1$ in $k^{th}$ place for each $k\in \mathbb{N}$ and $\phi=span\{e^{k}\}$, the set of all finitely non-zero sequences.\\

Let $A=(a_{nk})$ be an infinite matrix of complex numbers $a_{nk}$ and $x=(x_{k})\in \omega$, where $k,n\in\mathbb{N}$.
Then the sequence $Ax$ is called as the $A-$transform of $x$ defined by the usual matrix product.
Hence, we transform the sequence $x$ into the sequence $Ax=\{(Ax)_{n}\}$ where
\begin{eqnarray}\label{equ1}
(Ax)_{n}=\sum_{k}a_{nk}x_{k}
\end{eqnarray}
for each $n\in\mathbb{N}$, provided the series on the right hand side of (\ref{equ1}) converges for each $n\in\mathbb{N}$. Let $X$ and $Y$ be two sequence spaces. If $Ax$ exists and is in $Y$ for every sequence $x=(x_{k})\in X$, then we say that
$A$ defines a matrix mapping from $X$ into $Y$, and we denote it by writing $A :X \rightarrow Y$ if and only if the series on
the right hand side of (\ref{equ1}) converges for each $n\in\mathbb{N}$ and every $x\in X$, and we have $Ax=\{(Ax)_{n}\}_{n\in \mathbb{N}}\in Y$
for all $x\in X$.  A sequence $x$ is said to be $A$-summable to $l$ if $Ax$ converges to $l$ which is called the $A$-limit of $x$. Let $X$ be a sequence space and $A$ be an infinite matrix. The sequence space
\begin{eqnarray}\label{eq0}
X_{A}=\{x=(x_{k})\in\omega:Ax\in X\}
\end{eqnarray}
is called the domain of $A$ in $X$ which is a sequence space.\\

We write $\mathcal{U}$ for the set of all sequences $u=(u_{k})$ such that $u_{k}\neq 0$ for all $k\in \mathbb{N}$. For $u\in \mathcal{U}$, let $1/u=(1/u_{k})$. Let $u,w\in \mathcal{U}$. Now, we define the \emph{generalized weighted mean} or \emph{factorable matrix} $G(u,w)=(g_{nk})$ by

\begin{eqnarray*}
g_{nk}=\left\{\begin{array}{ccl}
u_{n}w_{k}&, & (0\leq k\leq n)\\
0&, & (k>n)
\end{array}\right.
\end{eqnarray*}

for all $k,n\in\mathbb{N}$; where $u_{n}$ depends only on $n$ and $w_{k}$ only on $k$.\\

By $\mathcal{F}$, we will denote the collection of all finite subsets on $\mathbb{N}$.
For simplicity in notation, here and in what follows, the summation without limits runs from $1$ to $\infty$.
Also we use the convention that any term with negative subscript is equal to zero.\\

\section{New Integrated and Differentiated Spaces}

In this section, we will give new spaces defined by a weighted mean.\\

The concepts of integrated and differentiated sequence spaces was firstly used by Goes and Goes \cite{Goes} as
$\int X=\left\{x=(x_{k})\in \omega: (kx_{k}) \in X \right\}$ and $d(X)=\left\{x=(x_{k})\in \omega: (k^{-1}x_{k}) \in X \right\}$.
Malkowsky and Sava\c{s} \cite{MalSav} have defined the sequence space $Z=(u,v;X)$, which consists of all sequences whose $G(u,v)-$ transforms are in
$X\in \{\ell_{\infty},c,c_{0},\ell_{p}\}$, where $u,w\in \mathcal{U}$. Paranormed sequence spaces derived by weighted mean are studied in \cite{AB5}. Altay and Ba\c{s}ar \cite{AB6} constructed the new paranormed sequence spaces $\ell(u,v;p)$. \c{S}imsek et al. \cite{PKS2} have introduced a modular structure of the sequence spaces defined by Altay and Ba\c{s}ar \cite{AB6} and studied Kadec-Klee and uniform Opial properties of this sequence space on K\"{o}the sequence spaces. In \cite{PKS}, using the generalized weighted mean, new difference sequence spaces are defined. Kiri\c{s}ci \cite{kirisci2} have defined the almost sequence spaces with generalized weighted mean and in \cite{kirisci}, studied some properties of new almost sequence spaces derived by generalized weighted mean. Structural properties of the $bv$ space are studied by Cesaro mean, generalized weighted mean and Riesz mean, in \cite{kirisci4}. Following the Goes and Goes \cite{Goes}, Kiri\c{s}ci \cite{kirisci3} have studied the integrated and differentiated sequence spaces and defined the Riesz type integrated and differentiated sequence spaces, in \cite{kirisci5}.\\

We define the new matrices $\Gamma=(\gamma_{nk})$ and $\Sigma=(\sigma_{nk})$ by
\begin{eqnarray}\label{matr1}
\gamma_{nk}= \left\{ \begin{array}{ccl}
ku_{n}\left(w_{k}-w_{k+1}\right)&, & \quad (k<n)\\
nu_{n}w_{n}&, & \quad (n=k)\\
0&, & \quad (k>n)
\end{array} \right.
\end{eqnarray}

\begin{eqnarray}\label{matr2}
\sigma_{nk}= \left\{ \begin{array}{ccl}
\frac{1}{k}u_{n}\left(w_{k}-w_{k+1}\right)&, & \quad (k<n)\\
\frac{u_{n}w_{n}}{n}&, & \quad (n=k)\\
0&, & \quad (k>n)
\end{array} \right.
\end{eqnarray}
for all $k,n\in\mathbb{N}$.\\

Let $u,w\in \mathcal{U}$. The new integrated spaces defined by
\begin{eqnarray*}
\int bv(u,w)=\left\{x=(x_{k})\in \omega: \sum_{k=1}^{n}u_{n}w_{k}\Delta(kx_{k})<\infty\right\}
\end{eqnarray*}
and the new differentiated spaces defined by
\begin{eqnarray*}
d(bv(u,w))=\left\{x=(x_{k})\in \omega: \sum_{k=1}^{n}u_{n}w_{k}\Delta(k^{-1}x_{k})<\infty\right\}.
\end{eqnarray*}

Consider the notation (\ref{eq0}) and the matrices (\ref{matr1}), (\ref{matr2}). From here, we can re-define the spaces $\int bv(u,w)$ and $d(bv(u,w))$ by
\begin{eqnarray}\label{domain}
\left(\ell_{1}\right)_{\Gamma}=\int bv(u,w)
\end{eqnarray}
and
\begin{eqnarray}\label{domain1}
\left(\ell_{1}\right)_{\Sigma}=d(bv(u,w)).
\end{eqnarray}

Let $x=(x_{k})\in \int bv(u,w)$. The $\Gamma-$transform of a sequence $x=(x_{k})$ is defined by
\begin{eqnarray}\label{transform1}
y_{n}=\sum_{k=1}^{n-1}ku_{n}(w_{k}-w_{k+1})x_{k}+nu_{n}w_{n}x_{n}
\end{eqnarray}
where $\Gamma$ is defined by (\ref{matr1}).
Let $x=(x_{k})\in d(bv(u,w))$. The $\Sigma-$transform of a sequence $x=(x_{k})$ is defined by
\begin{eqnarray}\label{transform2}
y_{n}=\sum_{k=1}^{n-1}\frac{1}{k}u_{n}(w_{k}-w_{k+1})x_{k}+\frac{1}{n}u_{n}w_{n}x_{n}
\end{eqnarray}
where $\Sigma$ is defined by (\ref{matr2}).

\begin{Theorem}\label{isomorph}
The integrated and differentiated sequence spaces derived by weighted mean are norm isomorphic to the absolute summable sequence space.
\end{Theorem}

\begin{proof}
We must show that a linear bijection between the integrated sequence space derived by weighted mean  and the absolute summable sequence space exists. Consider the transformation $f_{\Gamma}$ defined, with the notation (\ref{transform1}), from $\int bv(u,w)$ to $\ell_{1}$ by $x \mapsto y=f_{\Gamma}x$. The linearity of $f_{\Gamma}$ is clear. Also, it is trivial that $x=\theta$ whenever $f_{\Gamma}x=\theta$ and therefore, $f_{\Gamma}$ is injective.\\

Let $y\in \ell_{1}$ and define the sequence $x=(x_{k})$ by
\begin{eqnarray*}
x_{k}=\sum_{j=1}^{k-1}\frac{1}{k}\frac{1}{u_{j}}\left(\frac{1}{w_{j}}-\frac{1}{w_{j+1}}\right)y_{j}+\frac{y_{k}}{k.u_{k}w_{k}}.
\end{eqnarray*}
Then
\begin{eqnarray*}
\|x\|_{\int bv(u,w)}=\sum_{k}\left|\sum_{j=1}^{k-1}ju_{k}\left(w_{j}-w_{j+1}\right)x_{j}+nu_{n}w_{n}x_{n}\right|=\sum_{k}|y_{k}|=\|y\|_{\ell_{1}}<\infty.
\end{eqnarray*}
Then, we have that $x\in \int bv(u,w)$. So, $f_{\Gamma}$ is surjective and norm preserving. Hence $f_{\Gamma}$ is a linear bijection. It shown us that the space $\int bv(u,w)$ is norm isomorphic to $\ell_{1}$.\\

As similar, using the notation (\ref{transform2}), we can define the transformation $f_{\sum}$ from $d(bv(u,w))$ and $\ell_{1}$ by $x \mapsto y=f_{\sum}x$. If we choose the sequence $x=(x_{k})$ by
\begin{eqnarray*}
x_{k}=\sum_{j=1}^{k-1}k\frac{1}{u_{j}}\left(\frac{1}{w_{j}}-\frac{1}{w_{j+1}}\right)y_{j}+\frac{k.y_{k}}{u_{k}w_{k}}
\end{eqnarray*}
while $y\in \ell_{1}$, then we obtain the space $d(bv(u,w))$ is norm isomorphic to $\ell_{1}$
with the norm $\|x\|_{d(bv(u,w))}$.
\end{proof}

Since $\int bv(u,w)=\left[\ell_{1}\right]_{\Gamma}$ and $d(bv(u,w))=\left[\ell_{1}\right]_{\Sigma}$ holds, $\ell_{1}$ is a $BK-$space with the norm $\|x\|_{\ell_{1}}$ and the matrices $\Gamma$and $\Sigma$ are triangle matrix, then Theorem 4.3.2 of Wilansky\cite{Wil84} gives the fact that the integrated and differentiated sequence spaces derived by weighted mean are $BK-$space. Therefore, there is no need for detailed proof of the following theorem.

\begin{Theorem}
The spaces $\int bv(u,w)$ and $d(bv(u,w))$ are $BK-$space with the norms $\|x\|_{\int bv(u,w)}=\|x\|_{\ell_{1}(\Gamma)}$ and
$\|x\|_{d(bv(u,w))}=\|x\|_{\ell_{1}(\Sigma)}$, respectively.
\end{Theorem}

\begin{Theorem}
The differentiated sequence space derived by weighted mean has $AK-$property.
\end{Theorem}

\begin{proof}
Let $x=(x_{k})\in d(bv(u,w))$ and $x^{[n]}=\{x_{1}, x_{2}, \cdots, x_{n},0,0,\cdots\}$.
Hence,
\begin{eqnarray*}
x-x^{[n]}=\{0,0,\cdots,0,x_{n+1}, x_{n+2},\cdots\} \Rightarrow \|x-x^{[n]}\|_{d(bv(u,w))}=\|0,0,\cdots,0,x_{n+1}, x_{n+2},\cdots\|
\end{eqnarray*}
and since $x\in d(bv(u,w))$,
\begin{eqnarray*}
\|x-x^{[n]}\|_{d(bv(u,w))}=\sum_{k\geq n+1}\left|\frac{1}{k}u_{n}(w_{k}-w_{k+1})x_{k}+\frac{1}{n}u_{n}w_{n}x_{n}\right|\rightarrow 0 ~\textrm{ as}~\ n\rightarrow \infty\\
\Rightarrow \lim_{n\rightarrow \infty}\|x-x^{[n]}\|_{d(bv(u,w))}=0 \Rightarrow x^{[n]}\rightarrow x ~\textrm{ as}~\ n\rightarrow \infty ~\textrm{ in}~\ d(bv(u,w)).
\end{eqnarray*}
Then the space $d(bv(u,w))$ has $AK-$property.
\end{proof}

Because of the isomorphisms $f_{\Gamma}$ and $f_{\sum}$, defined in the proof of Theorem \ref{isomorph}, are onto the inverse image of the basis $\{e^{(k)}\}_{k\in \mathbb{N}}$ of the space $\ell_{1}$ is the basis of the spaces $\int bv(u,w)$ and $d(bv(u,w))$. Therefore, we can give following theorems for Schauder basis of new sequence spaces :

\begin{Theorem}\label{basisweighted}
Define a sequence $s^{(k)}=\{s_{n}^{(k)}\}_{n\in\mathbb{N}}$ of elements of the space $\int bv(u,w)$ for every fixed $k\in \mathbb{N}$ by
\begin{eqnarray*}
s_{n}^{(k)}= \left\{ \begin{array}{ccl}
\frac{1}{n}\left(\frac{1}{u_{k}w_{k}}-\frac{1}{u_{k}u_{k+1}}\right)&, & \quad (1< k <n)\\
\frac{1}{nu_{n}w_{n}}&, & \quad (n=k)\\
0&, & \quad (k>n)
\end{array} \right.
\end{eqnarray*}
Therefore, the sequence $\{s^{(k)}\}_{k\in\mathbb{N}}$ is a basis for the space $\int bv(u,w)$ and any $x\in \int bv(u,w)$ has a unique representation of the form
\begin{eqnarray}\label{basis1}
x=\sum_{k}(\Gamma x)_{k}s^{(k)}
\end{eqnarray}
\end{Theorem}

\begin{proof}
Let $e^{(k)}$ be a sequence whose only non-zero term is a $1$ in $k^{th}$ place for each $k\in \mathbb{N}$.
We know that
\begin{eqnarray}\label{basis2}
\Gamma s^{(k)}(q)=e^{(k)}\in \ell_{1}
\end{eqnarray}
for all $k\in\mathbb{N}$. Then, we have $\{s^{(k)}(q)\}\subset \int bv(u,w)$.\\

We take $x\in \int bv(u,w)$. Then, we put,
\begin{eqnarray}\label{basis3}
x^{[m]}=\sum_{k=1}^{m}(\Gamma x)_{k}(q)s^{(k)}(q),
\end{eqnarray}
for every positive integer $m$. Then, we have
\begin{eqnarray*}
\Gamma x^{[m]}=\sum_{k=1}^{m}(\Gamma x)_{k}(q)\Gamma s^{(k)}(q)=\sum_{k=1}^{m}(\Gamma x)_{k}(q)e^{(k)}
\end{eqnarray*}
and
\begin{eqnarray*}
\left(\Gamma(x-x^{[m]})\right)_{i}=\left\{ \begin{array}{ccl}
0&, & \quad (1\leq i <m)\\
(\Gamma x)_{i}&, & \quad (i>m)
\end{array} \right.
\end{eqnarray*}
by applying $\Gamma$ to (\ref{basis3}) with (\ref{basis2}), for $i,m\in \mathbb{N}$.
For $\varepsilon >0$, there exists an integer $m_{0}$ such that
\begin{eqnarray*}
\left[\sum_{i=m}^{\infty}|(\Gamma x)_{i}|\right]<\varepsilon / 2
\end{eqnarray*}
for all $m\geq m_{0}$. Hence,
\begin{eqnarray*}
\|x-x^{[m]}\|_{\int bv(u,w)}=
\end{eqnarray*}
for all $m\geq m_{0}$. Therefore, $x\in \int bv(u,w)$ is represented as in (\ref{basis1}), as we desired.
\end{proof}

\begin{Theorem}\label{basisweighted2}
Define a sequence $t^{(k)}=\{t_{n}^{(k)}\}_{n\in\mathbb{N}}$ of elements of the space $d(bv(u,w))$ for every fixed $k\in \mathbb{N}$ by
\begin{eqnarray*}
t_{n}^{(k)}= \left\{ \begin{array}{ccl}
n\left(\frac{1}{u_{k}w_{k}}-\frac{1}{u_{k}u_{k+1}}\right)&, & \quad (1< k <n)\\
\frac{n}{u_{n}w_{n}}&, & \quad (n=k)\\
0&, & \quad (k>n)
\end{array} \right.
\end{eqnarray*}
Therefore, the sequence $\{t^{(k)}\}_{k\in\mathbb{N}}$ is a basis for the space $d(bv(u,w))$ and any $x\in d(bv(u,w))$ has a unique representation of the form
\begin{eqnarray*}
x=\sum_{k}(\Sigma x)_{k}t^{(k)}.\\
\end{eqnarray*}
\end{Theorem}

\begin{remark}\label{rem}
It is well known that every Banach space $X$ with a Schauder basis is separable.
\end{remark}

From Theorem \ref{basisweighted}, Theorem \ref{basisweighted2} and Remark,
we can give following corollary:\\

\begin{Corollary}
The spaces $\int bv(u,w)$ and $d(bv(u,w))$ are separable.
\end{Corollary}

\section{Dual Spaces}

If $X,Y\subset \omega$ and $z$ any sequence, we can write $z^{-1}*X=\{x=(x_{k})\in \omega: xz\in X\}$ and $M(X,Y)=\bigcap_{x\in X}x^{-1}*Y$.
If we choose $Y=cs, bs$, then we obtain the $\beta-, \gamma- $duals of $X$, respectively as
\begin{eqnarray*}
X^{\beta}&=&M(X,cs)=\{z=(z_{k})\in \omega:  zx=(z_{k}x_{k})\in cs ~\textrm{for all }~   x\in X\}\\
X^{\gamma}&=&M(X,bs)=\{z=(z_{k})\in \omega:  zx=(z_{k}x_{k})\in bs ~\textrm{for all }~   x\in X\}.
\end{eqnarray*}

Let $A=(a_{nk})$ be an infinite matrix. Now we give some conditions:

\begin{eqnarray}\label{deq1}
&&\sup_{k,n \in \mathbb{N}}|a_{nk}| < \infty,\\ \label{deq2}
&&\lim_{n \rightarrow \infty} a_{nk} = \alpha_{k}~ \textrm{ for each }~k\in\mathbb{N},\\ \label{deq3}
&&\sup_{k \in \mathbb{N}}\sum_{n}|a_{nk}| < \infty \\ \label{deq4}
&&\sup_{k,m \in \mathbb{N}}\left|\sum_{n=0}^{m}a_{nk}\right| < \infty, \\ \label{deq5}
&&\sum_{n}a_{nk} ~ \textrm{ convergent for each }~k\in\mathbb{N} \\ \label{deq6}
&&\sum_{n}a_{nk}=0 ~ \textrm{ for each }~ k\in\mathbb{N}
\end{eqnarray}

\begin{Lemma}\label{duallem1}
For the characterization of the class $(X:Y)$ with
$X=\{\ell_{1}\}$ and $Y=\{\ell_{\infty}, c,\ell_{1}\}$, we can give the necessary and sufficient
conditions from Table 1, where

\begin{center}
\begin{tabular}{|c | c | c |c | c | c |}
\hline
\textbf{1.} (\ref{deq1}) & \textbf{2.} (\ref{deq1}), (\ref{deq2})  & \textbf{3.} (\ref{deq3}) & \textbf{4.} (\ref{deq4}) & \textbf{5.} (\ref{deq4}), (\ref{deq5}) & \textbf{6.} (\ref{deq4}), (\ref{deq6}) \\
\hline
\end{tabular}
\end{center}

\end{Lemma}

\begin{center}
\begin{tabular}{|c | c c c c c c |}
\hline
To $\rightarrow$ & $\ell_{\infty}$ & $c$ & $\ell_{1}$ & $bs$ & $cs$ & $c_{0}s$\\ \hline
From $\downarrow$ &  &  & & & &\\ \hline
$\ell_{1}$ & \textbf{1.} & \textbf{2.} & \textbf{3.} & \textbf{4.} & \textbf{5.} & \textbf{6.}\\
\hline
\end{tabular}

\vspace{0.1cm}Table 1\\

\end{center}

\begin{Theorem}
We define the matrix $E=(e_{nk})$ as
\begin{eqnarray}\label{alphamtrx}
e_{nk}= \left\{ \begin{array}{ccl}
\frac{1}{n}\frac{1}{u_{k}}\left(\frac{1}{w_{k}}-\frac{1}{w_{k+1}}\right)a_{n}&, & \quad (1\leq k<n)\\
\frac{a_{n}}{nu_{n}w_{n}}&, & \quad (n=k)\\
0&, & \quad (k>n)
\end{array} \right.
\end{eqnarray}
for all $k,n\in \mathbb{N}$, where $u,w\in \mathcal{U}, a=(a_{k})\in\omega$.
The $\alpha-$dual of the space $\int bv(u,w)$ is the set
\begin{eqnarray*}
d_{1}=\left\{a=(a_{k})\in \omega: \sup_{N\in \mathcal{F}}\sum_{k}\left|\sum_{n\in N}e_{nk}\right|<\infty \right\}
\end{eqnarray*}
\end{Theorem}

\begin{proof}
Let $a=(a_{k})\in \omega$. We can easily derive that with the notation (\ref{transform1}) that
\begin{eqnarray}\label{alphaequ}
a_{n}x_{n}=\sum_{k=1}^{n}\frac{1}{n}\frac{1}{u_{k}}\left(\frac{1}{w_{k}}-\frac{1}{w_{k+1}}\right)a_{n}y_{k}+\frac{a_{n}y_{n}}{nu_{n}w_{n}}=\sum_{k=1}^{n}e_{nk}y_{k}=(Ey)_{n}
\end{eqnarray}
for all $k,n\in \mathbb{N}$, where $E=(e_{nk})$ is defined by (\ref{alphamtrx}). It follows from (\ref{alphaequ}) that $ax=(a_{n}x_{n})\in\ell_{1}$ whenever $x\in \int bv(u,w)$ if and only if $Ey\in \ell_{1}$ whenever $y\in \ell_{1}$. We obtain that $a\in[\int bv(u,w)]^{\alpha}$ whenever $x\in \int bv(u,w)$ if and only if $E\in (\ell_{1}:\ell_{1})$. Therefore, we get by Lemma \ref{duallem1} with $E$ instead of $A$ that $a\in\left[\int bv(u,w)\right]^{\alpha}$ if and only if $\sup_{k\in \mathbb{N}}\sum_{n}\left|e_{nk}\right|<\infty$.
This gives us the result that $\left[\int bv(u,w)\right]^{\alpha}=d_{1}$.
\end{proof}

\begin{Theorem}
We define the matrix $F=(f_{nk})$ as
\begin{eqnarray}\label{alphamtrx1}
f_{nk}= \left\{ \begin{array}{ccl}
n\frac{1}{u_{k}}\left(\frac{1}{w_{k}}-\frac{1}{w_{k+1}}\right)a_{n}&, & \quad (1\leq k<n)\\
\frac{na_{n}}{u_{n}w_{n}}&, & \quad (n=k)\\
0&, & \quad (k>n)
\end{array} \right.
\end{eqnarray}
for all $k,n\in \mathbb{N}$, where $u,w\in \mathcal{U}, a=(a_{k})\in\omega$.
The $\alpha-$dual of the space $d(bv(u,w))$ is the set
\begin{eqnarray*}
d_{2}=\left\{a=(a_{k})\in \omega: \sup_{N\in \mathcal{F}}\sum_{k}\left|\sum_{n\in N}f_{nk}\right|<\infty \right\}
\end{eqnarray*}
\end{Theorem}

\begin{Theorem}
Let $u,w\in \mathcal{U}$ for all $k,n\in \mathbb{N}$. Then the $\beta-$dual of the space $\int bv(u,w)$ is $d_{3}\cap cs$, where
\begin{eqnarray*}
d_{3}=\left\{a=(a_{k})\in \omega: \sum_{k=1}^{n}\left|\frac{1}{k}\frac{a_{k}}{u_{k}w_{k}}+\left(\frac{1}{u_{k}w_{k}}-\frac{1}{u_{k}w_{k+1}}\right)\sum_{j=k+1}^{n}\frac{1}{j}a_{j}\right|<\infty \right\}
\end{eqnarray*}
\end{Theorem}

\begin{proof}
Consider the equation
\begin{eqnarray}\label{betaequ}
\sum_{k=1}^{n}a_{k}x_{k}&=&\sum_{k=1}^{n}a_{k}\left[\sum_{j=1}^{k-1}\frac{1}{k}\left(\frac{1}{u_{j}w_{j}}-\frac{1}{u_{j}w_{j+1}}\right)y_{j}+\frac{y_{k}}{ku_{k}w_{k}}\right]
\end{eqnarray}
\begin{eqnarray*}
=\sum_{k=1}^{n}\left|\frac{1}{k}\frac{a_{k}y_{k}}{u_{k}w_{k}}+\left(\frac{1}{u_{k}w_{k}}-\frac{1}{u_{k}w_{k+1}}\right)y_{k}\sum_{j=k+1}^{n}\frac{1}{j}a_{j}\right|=(H_{n}y)
\end{eqnarray*}
for all $n\in \mathbb{N}$, where the matrix $H=(h_{nk})$ is defined by
\begin{eqnarray}\label{betamtrx}
h_{nk}= \left\{ \begin{array}{ccl}
\left(\frac{1}{u_{k}w_{k}}-\frac{1}{u_{k}w_{k+1}}\right)\sum_{j=k+1}^{n}\frac{1}{j}a_{j}&, & \quad (k>n\\
\frac{1}{n}\frac{a_{n}}{u_{n}w_{n}}&, & \quad (n=k)\\
0&, & \quad (k<n)
\end{array} \right.
\end{eqnarray}
for all $k,n\in \mathbb{N}$. Therefore, we deduce from Lemma \ref{duallem1} with (\ref{betaequ})
that $ax=(a_{n}x_{n})\in cs$ whenever $x\in \int bv(u,w)$ if and only if $Hy\in c$ whenever $y\in \ell_{1}$. From (\ref{deq1}) and (\ref{deq2}), we have
\begin{eqnarray*}
\lim_{n}h_{nk}=\alpha_{k}  ~ \textrm{ and }~ \sup_{k}\sum_{n}\left|h_{nk}\right|<\infty
\end{eqnarray*}
which shows that $\left[\int bv(u,w)\right]^{\beta}=d_{3}\cap cs$.
\end{proof}

\begin{Theorem}
$\left[\int bv(u,w)\right]^{\gamma}=d_{3}$.
\end{Theorem}

\begin{proof}
We obtain from Lemma \ref{duallem1} with (\ref{betaequ}) that $ax=(a_{n}x_{n})\in bs$ whenever $x\in \int bv(u,w)$ if and only if $Hy\in \ell_{\infty}$ whenever $y\in \ell_{1}$. Then, we see from (\ref{deq1}) that $\left[\int bv(u,w)\right]^{\gamma}=d_{3}$.\\
\end{proof}

\begin{Theorem}
The $\beta-$dual of the space $d(bv(u,w))$ is $d_{4}\cap cs$, where
\begin{eqnarray*}
d_{4}=\left\{a=(a_{k})\in \omega: \sum_{k=1}^{n}\left|\frac{k.a_{k}}{u_{k}w_{k}}+\left(\frac{1}{u_{k}w_{k}}-\frac{1}{u_{k}w_{k+1}}\right)\sum_{j=k+1}^{n}j.a_{j}\right|<\infty \right\}
\end{eqnarray*}
\end{Theorem}

\begin{Theorem}
$\left[d(bv(u,w))\right]^{\gamma}=d_{4}$.
\end{Theorem}

\section{Matrix transformations}

We shall write for brevity that
\begin{eqnarray*}
&&\overline{a}_{nk}=\sum_{k=1}^{n}\left|\frac{1}{k}\frac{a_{nk}}{u_{k}w_{k}}+\left(\frac{1}{u_{k}w_{k}}-\frac{1}{u_{k}w_{k+1}}\right)\sum_{j=k+1}^{n}\frac{1}{j}a_{nj}\right|,\\
&&\widetilde{a}_{nk}=\sum_{k=1}^{n}\left|\frac{k.a_{nk}}{u_{k}w_{k}}+\left(\frac{1}{u_{k}w_{k}}-\frac{1}{u_{k}w_{k+1}}\right)\sum_{j=k+1}^{n}j.a_{nj}\right|,\\
&&\overline{b}_{nk}=\sum_{j=1}^{n-1}j.u_{n}(w_{j}-w_{j+1})a_{jk}+n.u_{n}w_{n}a_{nk},\\
&&\widetilde{b}_{nk}=\sum_{j=1}^{n-1}\frac{1}{j}.u_{n}(w_{j}-w_{j+1})a_{jk}+\frac{1}{n}.u_{n}w_{n}a_{nk}
\end{eqnarray*}
for all $k,n\in\mathbb{N}$.

\begin{Theorem}\label{mtrxtr1}
Suppose that the entries of the infinite matrices $A=(a_{nk})$ and $B=(b_{nk})$ are connected with the relation
\begin{eqnarray}\label{trnsf2}
a_{nk}=\sum_{j=k}^{\infty} j.(w_{k}-w_{k+1})u_{j}b_{nj} \quad \quad ~\textrm{or }~ b_{nk}=\overline{a}_{nk}
\end{eqnarray}
for all $k,n\in\mathbb{N}$ and $Y$ be any given sequence space. Then $A\in (\int bv(u,w):Y)$ if and only if $\{a_{nk}\}_{k\in\mathbb{N}}\in \{\int bv(u,w)\}^{\beta}$ for all $n\in\mathbb{N}$ and $B\in (\ell_{1}:Y)$.\\
\end{Theorem}

\begin{proof}
Let $Y$ be any given sequence. Suppose that (\ref{trnsf2}) holds between the infinite matrices $A=(a_{nk})$ and $B=(b_{nk})$, and take into account that the spaces $\int bv(u,w)$ and $\ell_{1}$ are linearly isomorphic.

Let $A\in (\int bv(u,w):Y)$ and take any $y=(y_{k})\in \ell_{1}$. Then $B\Gamma$ exists and $\{a_{nk}\}_{k\in\mathbb{N}}\in \{\int bv(u,w)\}^{\beta}$ which yields that (\ref{trnsf2}) is necessary and $\{b_{nk}\}_{k\in\mathbb{N}}\in \ell_{1}^{\beta}$ for each $n\in \mathbb{N}$. Hence, $By$ exists for each $y\in \ell_{1}$ and thus by letting $m\rightarrow \infty$ in the equality
\begin{eqnarray*}
\sum_{k=1}^{m}a_{nk}x_{k}=\sum_{k=1}^{m}\left|\frac{1}{k}\frac{a_{nk}y_{k}}{u_{k}w_{k}}+\left(\frac{1}{u_{k}w_{k}}-\frac{1}{u_{k}w_{k+1}}\right)y_{k}\sum_{j=k+1}^{m}\frac{1}{j}a_{nj}\right| ~ \textrm{ for all }~m,n\in\mathbb{N}
\end{eqnarray*}
we obtain that $Ax=By$ which leads us to the consequence $B\in (\ell_{1}:Y)$.

Conversely, let $\{a_{nk}\}_{k\in\mathbb{N}}\in \{\int bv(u,w)\}^{\beta}$ for each $n\in \mathbb{N}$ and $B\in (\ell_{1}:Y)$, and take any $x=(x_{k})\in \int bv(u,w)$. Then, $Ax$ exists. Therefore, we obtain from the equality
\begin{eqnarray*}
\sum_{k=1}^{m}b_{nk}y_{k}=\sum_{k=1}^{m}a_{nk}x_{k} ~ \textrm{ for all }~m,n\in \mathbb{N}
\end{eqnarray*}
as $m\rightarrow\infty$ the result that $By=Ax$ and this shows that $A\in (\int bv(u,w):Y)$. This completes the proof.
\end{proof}

\begin{Theorem}\label{mtrxtr1d}
Suppose that the entries of the infinite matrices $A=(a_{nk})$ and $C=(c_{nk})$ are connected with the relation $c_{nk}=\overline{b}_{nk}$
for all $k,n\in\mathbb{N}$ and $Y$ be any given sequence space. Then, $A\in (Y:\int bv(u,w))$ if and only if $C\in (Y:\ell_{1})$.\\
\end{Theorem}

\begin{proof}
Let $z=(z_{k})\in Y$ and consider the following equality:
\begin{eqnarray}\label{trnsf1}
\sum_{k=1}^{m}c_{nk}z_{k}=\sum_{j=1}^{n}j.u_{n}(w_{j}-w_{j+1})\left(\sum_{k=1}^{m}a_{jk}z_{k}\right)
\end{eqnarray}
for all $m,n\in \mathbb{N}$. Equation (\ref{trnsf1}) yields as $m\rightarrow \infty$ the result that
$(Cz)_{n}=\{\Gamma(Az)\}_{n}$. Therefore, one can immediately observe from this that $Az\in \int bv(u,w)$
whenever $z\in Y$ if and only if $Cz\in \ell_{1}$ whenever $z\in Y$.
\end{proof}

\begin{Theorem}\label{mtrxtr2}
Suppose that the entries of the infinite matrices $A=(a_{nk})$ and $D=(d_{nk})$ are connected with the relation
\begin{eqnarray*}
a_{nk}=\sum_{j=k}^{\infty} \frac{1}{j}.(w_{k}-w_{k+1})u_{j}d_{nj} \quad \quad ~\textrm{or }~ d_{nk}=\widetilde{a}_{nk}
\end{eqnarray*}
for all $k,n\in\mathbb{N}$ and $Y$ be any given sequence space. Then $A\in (d(bv(u,w)):Y)$ if and only if $\{a_{nk}\}_{k\in\mathbb{N}}\in \{d(bv(u,w))\}^{\beta}$ for all $n\in\mathbb{N}$ and $D\in (\ell_{1}:Y)$.\\
\end{Theorem}

\begin{Theorem}\label{mtrxtr2d}
Suppose that the entries of the infinite matrices $A=(a_{nk})$ and $E=(e_{nk})$ are connected with the relation $e_{nk}=\widetilde{b}_{nk}$
for all $k,n\in\mathbb{N}$ and $Y$ be any given sequence space. Then, $A\in (Y:d(bv(u,w)))$ if and only if $E\in (Y:\ell_{1})$.\\
\end{Theorem}

\section{Examples}

\begin{Example}\label{exmpE}
The Euler sequence space $e_{\infty}^{r}$ is defined by $e_{\infty}^{r}=\{x\in \omega: \sup_{n\in\mathbb{N}}|\sum_{k=0}^{n}\binom{n}{k}(1-r)^{n-k}r^{k}x_{k}|<\infty\}$ (\cite{BF2}).
We consider the infinite matrix $A=(a_{nk})$ and define the matrix $F=(f_{nk})$ by
\begin{eqnarray*}
f_{nk}=\sum_{j=0}^{n}\binom{n}{j}(1-r)^{n-j}r^{j}a_{jk}  \quad \quad (k,n\in \mathbb{N}).
\end{eqnarray*}
If we want to get necessary and sufficient conditions for the class $(\int bv (u,w): e_{\infty}^{r})$ in Theorem \ref{mtrxtr1},
then, we replace the entries of the matrix $A$ by those of the matrix $F$.
\end{Example}

\begin{Example}\label{exmpR}
Let $T_{n}=\sum_{k=0}^{n}t_{k}$ and $A=(a_{nk})$ be an infinite matrix. We define the matrix $H=(h_{nk})$ by
\begin{eqnarray*}
h_{nk}=\frac{1}{T_{n}}\sum_{j=0}^{n}t_{j}a_{jk}  \quad \quad (k,n\in \mathbb{N}).
\end{eqnarray*}
Then, the necessary and sufficient conditions in order for $A$ belongs to the class $(\int bv(u,w):r_{\infty}^{t})$
are obtained from in Theorem \ref{mtrxtr1} by replacing the entries of the matrix $A$ by those of the matrix $H$;
 where $r_{\infty}^{t}$ is the space of all sequences whose $R^{t}-$transforms is in the space $\ell_{\infty}$ \cite{malk}.
\end{Example}

\begin{Example}
In the space $r_{\infty}^{t}$, if we take $t=e$, then, this space become to the Ces\`{a}ro sequence space of non-absolute type $X_{\infty}$ \cite{NgLee}.
As a special case, Example \ref{exmpR} includes the characterization of class $(\int bv(u,w):r_{\infty}^{t})$.
\end{Example}

\begin{Example}
The Taylor sequence space $t_{\infty}^{r}$ is defined by $t_{\infty}^{r}=\{x\in \omega: \sup_{n\in\mathbb{N}}|\sum_{k=n}^{\infty}\binom{k}{n}(1-r)^{n+1}r^{k-n}x_{k}|<\infty\}$ (\cite{kiris2}).
We consider the infinite matrix $A=(a_{nk})$ and define the matrix $P=(p_{nk})$ by
\begin{eqnarray*}
p_{nk}=\sum_{k=n}^{\infty}\binom{k}{n}(1-r)^{n+1}r^{k-n}a_{jk}  \quad \quad (k,n\in \mathbb{N}).
\end{eqnarray*}
If we want to get necessary and sufficient conditions for the class $(\int bv(u,w): t_{\infty}^{r})$ in Theorem \ref{mtrxtr1},
then, we replace the entries of the matrix $A$ by those of the matrix $P$.
\end{Example}

Similar to above examples, we can give necessary and sufficient conditions for the class $(d(bv(u,w)): Y)$ in Theorem \ref{mtrxtr2},
where $Y\in \{e_{\infty}^{r}, r_{\infty}^{t}, X_{\infty}, t_{\infty}^{r}\}$.\\

If we take the spaces $\ell_{\infty}$, $c$, $c_{0}$, $bs$, $cs$ and $c_{0}s$ instead of $X$ in Theorem \ref{mtrxtr2}, or $Y$ in Theorem \ref{mtrxtr1}
we can write the following examples. Firstly, we give some conditions and following lemmas:

\begin{eqnarray}\label{deq0}
&&\sup_{N,K\in \mathcal{F}}\left|\sum_{n\in N}\sum_{k\in K}a_{nk}\right|<\infty, \\\label{deq7}
&&\lim_{k}a_{nk}=0 ~ \textrm{ for each }~ n\in\mathbb{N}, \\ \label{deq8}
&&\sup_{N,K\in \mathcal{F}}\left|\sum_{n\in N}\sum_{k\in K}(a_{nk}-a_{n,k+1})\right|<\infty, \\ \label{deq9}
&&\sup_{N,K\in \mathcal{F}}\left|\sum_{n\in N}\sum_{k\in K}(a_{nk}-a_{n,k-1})\right|<\infty
\end{eqnarray}

\begin{Lemma}\label{lemmtr1}
Consider that the $X\in\{\ell_{\infty}, c, c_{0}, bs, cs, c_{0}s\}$ and $Y\in \{\ell_{1}\}$.
The necessary and sufficient conditions for $A\in (X:Y)$ can be read the following, from Table 2:

\begin{center}
\begin{tabular}{|c | c | c |c | }
\hline
\textbf{7.} (\ref{deq0}) & \textbf{8.} (\ref{deq7}), (\ref{deq8})  & \textbf{9.} (\ref{deq9}) & \textbf{10.} (\ref{deq8}) \\
\hline
\end{tabular}
\end{center}

\end{Lemma}

\begin{center}
\begin{tabular}{|c | c c  c c c c |}
\hline
From $\rightarrow$ & $\ell_{\infty}$ & $c$ & $c_{0}$ & $bs$ & $cs$, & $c_{0}s$\\ \hline
To $\downarrow$ &    & & & & &\\ \hline
$\ell_{1}$ & \textbf{7.} & \textbf{7.} & \textbf{7.} & \textbf{8.} & \textbf{9.} & \textbf{10.}\\
\hline
\end{tabular}

\vspace{0.1cm}Table 2\\
\end{center}

\begin{Example}
We choose $X\in \{\int bv(u,w)\}$ and $Y\in \{\ell_{\infty}, c, c_{0}\}$.
The necessary and sufficient conditions for $A\in (X:Y)$ can be taken from the Table 3:
\end{Example}
\begin{itemize}
  \item[\textbf{1a.}] (\ref{deq1}) holds with $\overline{a}_{nk}$ instead of ${a}_{nk}$.
\item[\textbf{2a.}] (\ref{deq1}) and  (\ref{deq2}) hold with $\overline{a}_{nk}$ instead of ${a}_{nk}$.
\item[\textbf{3a.}] (\ref{deq1}) and  (\ref{deq2}) hold with $\alpha_{k}=0$ as $\overline{a}_{nk}$ instead of $a_{nk}$.
\item[\textbf{4a.}] (\ref{deq4}) holds with $\overline{a}_{nk}$ instead of ${a}_{nk}$
\item[\textbf{5a.}] (\ref{deq4}), (\ref{deq5}) hold with $\overline{a}_{nk}$ instead of ${a}_{nk}$.
\item[\textbf{6a.}] (\ref{deq4}), (\ref{deq6}) hold with $\overline{a}_{nk}$ instead of ${a}_{nk}$.
\end{itemize}

\begin{center}
\begin{tabular}{|c | c c c c c c|}
\hline
To $\rightarrow$ & $\ell_{\infty}$ & $c$ & $c_{0}$ & $bs$ & $cs$ & $c_{0}s$ \\ \hline
From $\downarrow$ &  &  & & & &\\ \hline
$\int bv(u,w)$ & \textbf{1a.} & \textbf{2a.} & \textbf{3a.} & \textbf{4a.} & \textbf{5a.} & \textbf{6a.}\\
\hline
\end{tabular}

\vspace{0.1cm}Table 3\\
\end{center}

\begin{Example}
We choose $X\in \{d(bv(u,w))\}$ and $Y\in \{\ell_{\infty}, c, c_{0}, bs, cs, c_{0}s\}$.
The necessary and sufficient conditions for $A\in (X:Y)$ can be taken from the Table 4:
\end{Example}
\begin{itemize}
  \item[\textbf{1b.}] (\ref{deq1}) holds with $\widetilde{a}_{nk}$ instead of ${a}_{nk}$.
\item[\textbf{2b.}] (\ref{deq1}) and  (\ref{deq2}) hold with $\widetilde{a}_{nk}$ instead of ${a}_{nk}$.
\item[\textbf{3b.}] (\ref{deq1}) and  (\ref{deq2}) hold with $\alpha_{k}=0$ as $\widetilde{a}_{nk}$ instead of $a_{nk}$.
\item[\textbf{4b.}] (\ref{deq4}) holds with $\widetilde{a}_{nk}$ instead of ${a}_{nk}$
\item[\textbf{5b.}] (\ref{deq4}), (\ref{deq5}) hold with $\widetilde{a}_{nk}$ instead of ${a}_{nk}$.
\item[\textbf{6b.}] (\ref{deq4}), (\ref{deq6}) hold with $\widetilde{a}_{nk}$ instead of ${a}_{nk}$.
\end{itemize}

\begin{center}
\begin{tabular}{|c | c c c c c c |}
\hline
To $\rightarrow$ & $\ell_{\infty}$ & $c$ & $c_{0}$ & $bs$ & $cs$ & $c_{0}s$\\ \hline
From $\downarrow$ &  &  & & & &\\ \hline
$d(bv(u,w))$ & \textbf{1b.} & \textbf{2b.} & \textbf{3b.} & \textbf{4b.} & \textbf{5b.} & \textbf{6b.}\\
\hline
\end{tabular}

\vspace{0.1cm}Table 4\\
\end{center}

Using the Lemma \ref{lemmtr1}, we can give the Table 5 for $X\in \{\ell_{\infty}, c, c_{0}, bs, cs, c_{0}s\}$ and $Y\in \{\int bv(u,w)\}$
and Table 6 for $X\in \{\ell_{\infty}, c, c_{0}, bs, cs, c_{0}s\}$ and $Y\in \{d(bv(u,w))\}$ as follows:

\begin{itemize}
  \item[\textbf{7a.}] (\ref{deq0}) hold with $\overline{b}_{nk}$ instead of ${a}_{nk}$.
\item[\textbf{8a.}] (\ref{deq7}) and (\ref{deq8}) hold with $\overline{b}_{nk}$ instead of ${a}_{nk}$.
\item[\textbf{9a.}] (\ref{deq9}) holds with $\overline{b}_{nk}$ instead of ${a}_{nk}$.
\item[\textbf{10a.}] (\ref{deq8}) holds with $\overline{b}_{nk}$ instead of ${a}_{nk}$.
  \item[\textbf{7b.}] (\ref{deq0}) hold with $\widetilde{b}_{nk}$ instead of ${a}_{nk}$.
\item[\textbf{8b.}] (\ref{deq7}) and (\ref{deq8}) hold with $\widetilde{b}_{nk}$ instead of ${a}_{nk}$.
\item[\textbf{9b.}] (\ref{deq9}) holds with $\widetilde{b}_{nk}$ instead of ${a}_{nk}$.
\item[\textbf{10b.}] (\ref{deq8}) holds with $\widetilde{b}_{nk}$ instead of ${a}_{nk}$.
\end{itemize}

\begin{center}
\begin{tabular}{|c | c c  c c c c |}
\hline
From $\rightarrow$ & $\ell_{\infty}$ & $c$ & $c_{0}$ & $bs$ & $cs$, & $c_{0}s$\\ \hline
To $\downarrow$ &    & & & & &\\ \hline
$\int bv(u,w)$ & \textbf{7a.} & \textbf{7a.} & \textbf{7a.} & \textbf{8a.} & \textbf{9a.} & \textbf{10a.}\\
\hline
\end{tabular}

\vspace{0.1cm}Table 5\\
\end{center}

\begin{center}
\begin{tabular}{|c | c c  c c c c |}
\hline
From $\rightarrow$ & $\ell_{\infty}$ & $c$ & $c_{0}$ & $bs$ & $cs$, & $c_{0}s$\\ \hline
To $\downarrow$ &    & & & & &\\ \hline
$d(bv(u,w))$ & \textbf{7b.} & \textbf{7b.} & \textbf{7b.} & \textbf{8b.} & \textbf{9b.} & \textbf{10b.}\\
\hline
\end{tabular}

\vspace{0.1cm}Table 6\\
\end{center}

\section{Conclusion}
The difference sequence spaces are given by K{\i}zmaz \cite{Kizmaz}. If we choose the absolute summable
sequence space and apply the difference operator to this space, we obtain the space of all sequences of bounded variation and denote
by $bv$. The space $bv_{p}$ consisting of all sequences whose differences are in the space $\ell_{p}$.
The space $bv_{p}$ was introduced by Ba\c{s}ar and Altay \cite{AB2}. More recently, the sequence spaces $bv$ and $bv_{p}$
are studied in \cite{AB2}, \cite{BAM}, \cite{imamiri}, \cite{JarMal}, \cite{kirisci4}, \cite{kirisci3}, \cite{kirisci5}, \cite{MalRakZiv}.\\

Integrated and differentiated sequence spaces are introduced by \cite{Goes}. Kiri\c{s}ci \cite{kirisci3} have studied some properties of these spaces and
defined the Riesz type integrated and differentiated sequence spaces \cite{kirisci5}. In this work, we define the new integrated and differentiated sequence
spaces. We also compute the dual spaces of these spaces which are help us in the characterization of matrix mappings. Therefore, we characterize the matrix classes.
In last section, we give some examples related to the matrix transformations
in the table form.

\section*{Competing Interests}

The author declares that they have no competing interests.

{\bf Received: Month xx, 20xx}


\begin{thebibliography}{99}

\bibitem{AB2} F. Ba\c{s}ar and B. Altay \textit{On the space of sequences of $p$-bounded variation and related matrix mappings},
Ukranian Math. J., \textbf{55(1)} (2003), 136--147.

\bibitem{BF2} B. Altay, F. Ba\c sar, \textit{On some Euler sequence spaces
of non-absolute type}, Ukrainian Math. J. \textbf{57}(1)(2005), 1--17.

\bibitem{AB5} B. Altay and F. Ba\c{s}ar, \textit{Some paranormed sequence spaces of non-absolute type derived by weighted mean},
J. Math. Analysis and Appl., \textbf{319} (2006), 494--508.

\bibitem{AB6} B. Altay and F. Ba\c{s}ar, \textit{Generalization of the sequence space $\ell(p)$ derived by weighted mean},
J. Math. Anal. Appl., \textbf{330} (2007), 174--185.


\bibitem{Basarkitap} F. Ba\c{s}ar, \textit{Summability Theory and its Applications}, Bentham Science Publishers, e-books, Monographs, (2011).

\bibitem{BAM} F. Ba\c{s}ar, B. Altay M. Mursaleen, \textit{Some generalizations of the sequence space $bv_{p}$ of $p-$bounded variation sequences},
Nonlinear Analysis, \textbf{68} (2008), 273--287.

\bibitem{Goes} G. Goes and S., Goes,  \textit{Sequences of bounded variation and sequences of Fourier coefficients I},
Math. Z., \textbf{118}(1970), 93--102.

\bibitem{imamiri} M. Imaninezhad and M. Miri,  \textit{The dual space of the sequence space $bv(p)$, $1\leq p < \infty$},
Acta Math. Univ. Comenianae, \textbf{79}(2010), 143--149.

\bibitem{JarMal} A.M. Jarrah and  E. Malkowsky,  \textit{The space $bv(p)$, its $\beta-$dual and matrix transformations},
Collect. Math., \textbf{55(2)}(2004), 151--162.

\bibitem{Kizmaz} H. K{\i}zmaz,  \textit{On certain sequence spaces},
Can. Math. Bull., \textbf{24(2)}(1981), 169--176.

\bibitem{kirisci2} M. Kiri\c{s}ci, \textit{Almost convergence and generalized weighted mean},
First International Conference on Analysis and Applied Mathematics, AIP Conf. Proc. \textbf{1470}, 191--194, (2012), http://dx.doi.org/10.1063/1.4747672.

\bibitem{kirisci} M. Kiri\c{s}ci, \textit{Almost convergence and generalized weighted mean II},
Journal of Inequalities and Applications 2014, 2014:93, doi: 10.1186/1029-242X-2014-93.

\bibitem{kirisci4} M. Kiri\c{s}ci, \textit{The space $bv$ and some applications},
Mathematica Aeterna, \textbf{4(3)}(2014), 207--223.

\bibitem{kirisci3} M. Kiri\c{s}ci, \textit{Integrated and differentiated sequence spaces},
J. Nonlinear Anal. Appl., \textbf{1}(2015), 2--16, doi:10.5899/2015/jnna-00266.

\bibitem{kirisci5} M. Kiri\c{s}ci, \textit{Riesz type integriated and differentiated sequence spaces},
Bull. Math. Anal. Appl. , \textbf{2}(2015), 14--27.

\bibitem{kiris2} M. Kiri\c s\c ci, \textit{On the Taylor sequence spaces of nonabsolute type which include the spaces $c_{0}$ and $c$},
J. Math. Anal., \textbf{6(2)}(2015), 22--35.

\bibitem{malk} Malkowsky,E.:
\textit{Recent results in the theory of matrix transformations in sequence spaces}, Mat. Ves. \textbf{49},
187--196, (1997).

\bibitem{MalRakZiv} E. Malkowsky, V. Rakocevic and S. Zivkovic,  \textit{Matrix transformations between the
 sequence space $bv^{p}$ and certain BK spaces},
Bull. Cl. Sci. Math. Nat. Sci. Math., \textbf{27}(2002), 331--46.

\bibitem{MalSav} E. Malkowsky and E. Sava\c{s},  \textit{Matrix transformations between sequence spaces of generalized weighted means},
Appl. Math. Comput., \textbf{147(2)}(2012), 333--345.

\bibitem{NgLee} Ng, P.-N., Lee, P. -Y.: \textit{Ces\`{a}ro sequence spaces of non-absolute type}, Comment. Math. Prace Mat. \textbf{20(2)}, 429--433, (1978).

\bibitem{PKS2} N.\c{S}im\c{s}ek, V. Karakaya and H. Polat \textit{On some geometrical properties of generalized modular spaces of Cesaro type
defined by weighted means},
J. Inequal. Appl., 2009, Article ID 932734 (2009).

\bibitem{PKS} H. Polat, V. Karakaya and N.\c{S}im\c{s}ek \textit{Difference sequence spaces derived by generalized weighted mean},
Appl. Math. Lett., \textbf{24(5)} (2011), 608--614.

\bibitem{Wil84} A. Wilansky, \textit{Summability through Functinal Analysis}, North Holland, New York, (1984).
\end{thebibliography}
\end{document}